%% file: musicology_name_1.tex
\documentclass{amsart}

\usepackage{epstopdf}
\usepackage[bookmarksnumbered,colorlinks,bookmarks,citecolor=blue,linkcolor=blue,urlcolor=blue,breaklinks,linktocpage]{hyperref}
\usepackage[all]{hypcap}
\usepackage{subfigure}
\usepackage{graphicx}

\usepackage[all,cmtip]{xy}
\usepackage{calc}
\usepackage{enumitem}
\usepackage{multirow}

\newcommand\fini{\hfill $\clubsuit$}

\DeclareMathOperator{\sym}{\mathbf{Sym}}
\let\x\delta
\DeclareMathOperator{\e}{\epsilon}

\theoremstyle{plain}
\newtheorem{theorem}{Theorem}[section]

\theoremstyle{definition}
\newtheorem{definition}[theorem]{Definition}

\theoremstyle{remark}

\begin{document}



\title[Musicological aspects of counterpoint theory]{Musicological, computational, and conceptual aspects of first-species counterpoint theory}

\author[J. S. Arias-Valero]{Juan Sebasti\'{a}n Arias-Valero}
\address{Departamento de Matem\'{a}ticas, Universidad Nacional Aut\'{o}noma de M\'{e}xico, Ciudad de M\'{e}xico, M\'{e}xico}
\email{jsariasv1@gmail.com}
\author[O. A. Agustín-Aquino]{Octavio Alberto Agustín-Aquino}
\address{Instituto de Física y Matemáticas, Universidad Tecnológica de la Mixteca, Huajuapan de León, México}
\email{octavioalberto@mixteco.utm.mx}
\author[E. Lluis-Puebla]{Emilio Lluis-Puebla}
\address{Departamento de Matem\'{a}ticas, Universidad Nacional Aut\'{o}noma de M\'{e}xico, Ciudad de M\'{e}xico, M\'{e}xico}
\email{lluisp@unam.mx}

\begin{abstract}
We re-create the essential results of a 1989 unpublished article by Mazzola and Muzzulini that contains the musicological aspects of a first-species counterpoint model. We include a summary of the mathematical counterpoint theory and several variations of the model that offer different perspectives on Mazzola's original principles.
\end{abstract}





\maketitle

\section{Introduction}\label{intro}

The original idea of this article was to communicate the musicological results from \cite{MazzMuzz}, which presents the results of a model for \textit{first-species counterpoint}, based on ring and module theory, but was not published. This kind of counterpoint is the simplest one and the didactic basis of \textit{Renaissance counterpoint}, as taught by \cite{Fux}. The model was introduced by \cite{Inicount} for $\mathbb{Z}_{12}$; a ring that can be used to model the algebraic behavior of the twelve intervals between tones in the \textit{chromatic scale} of Western musical tradition. Then, the model was re-exposed with some additional computational results by Hichert in \cite[Part~VII]{MazzolaTopos}. Further generalizations to the case when the rings are of the form $\mathbb{Z}_{n}$ were considered in \cite{OAAAthesis,Junod}, and then included in a collaborative compendium of mathematical counterpoint theory and its computational aspects \cite{Octavio}. The motivation for such a generalization to $\mathbb{Z}_n$ was the existence of \textit{microtonal scales} with more than twelve tones, which have been used for making real music \cite{octaviotod}. 

However, the understanding of \cite{MazzMuzz} is unavoidably linked to that of the very model, so in this paper we expose a summary of it and the musicological, computer-aided, study of its results. We do not reduce this study to a mere exposition of the previous results, but propose several variations of the model with the due justifications. The intention of these variations is not to destroy the original model, but to deepen into some of its basic principles. 

Regarding other mathematical models of counterpoint, Mazzola's school has an important counterpart in D. Tymoczko. Tymoczko's model \cite[Appendix]{Tympolem} is based on orbifolds and is oriented towards voice leading. Moreover, it works by means of a geometric reading of the usual counterpoint rules, in contrast to the predictive character of Mazzola's model, which follows the principle that these rules obey mathematical relations based on musical symmetries. There has been a polemic around Mazzola's and Tymoczko's models. Tymoczko's initial critique can be found in \cite{Tympolem} and a response by Mazzola and the second author can be found in \cite{Mazzpolem}. Two of Tymoczko's concerns just were that the musicological study \cite{MazzMuzz} was not available and that there were some progressions forbidden by the model but good to Fux. Likely, this article can clarify those and other possible concerns. 

We organize this article as follows. First, in Section~\ref{strsty} we codify Fux's \textit{strict style} following two standard sources \cite{Fux,Jep}. This is a \textit{descriptive model}, which helps to make a mathematical taxonomy of progressions into inadmissible, bad, and good ones, but hardly explains anything on the nature of the rules. Then, based on \cite[Chapitre~III]{Beau}, in Sections \ref{modprin} and \ref{mod}, we expose Mazzola's model both in informal and formal synthetic ways. The quantitative computational results are also included. In Section~\ref{redstrsty}, we reduce, modulo $12$, the strict style to obtain the \textit{reduced strict style}. This procedure helps to compare the original phenomenon to the model, whose base is the integers modulo $12$. Here, we note that there are some new \textit{ambiguous} progressions that probably do not deserve the defined names inadmissible, bad, or good. In Section~\ref{comp}, we compare the reduced strict style to the model. Based on the conceptual and structural understanding of the model recorded in the previous sections, in Section~\ref{natvar} we propose an alternative model that replaces the dual numbers product (infinitesimals one) by an integers modulo $12$ one. The alternative to the dual numbers ring is also isomorphic to the product $\mathbb{Z}_{12}\times \mathbb{Z}_{12}$, and therefore could be suitable for further categorical generalizations. The quantitative results of the alternative model are very similar to those of the original one, and slightly improve them regarding the predictions of inadmissible or bad progressions in the reduced strict style. Also, in Sections~\ref{2var} and \ref{3var}, we propose two variations of the model that deepen into the principle of local characterization of deformed consonances/dissonances. These variations become the same and, again,  they slightly improve the prediction of inadmissible or bad progressions. Finally, we provide conclusions from our study and point out some directions for further research. 

\section{The strict style}\label{strsty}

\begin{figure}
\centering
\def\svgscale{0.8}
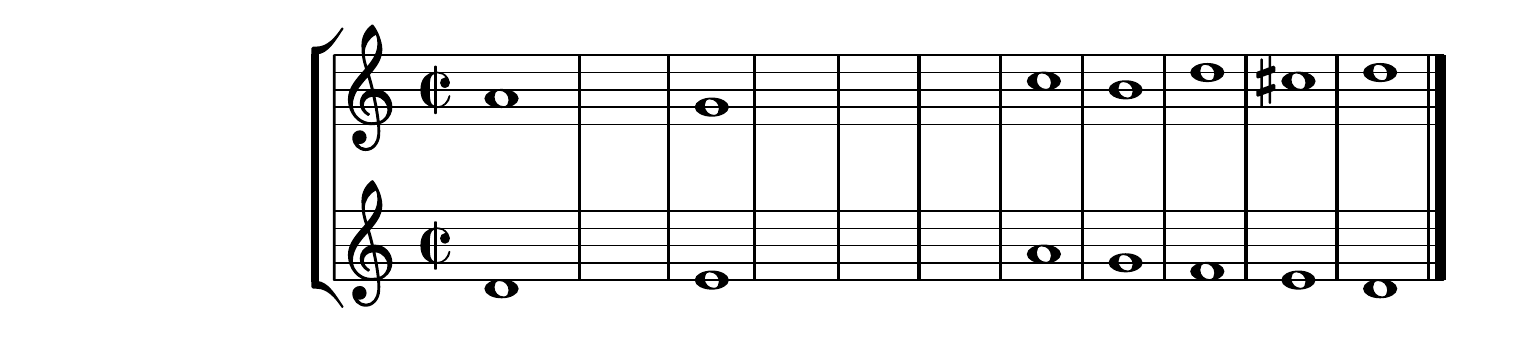
\caption{Some features of a first-species counterpoint example given by \cite[p.~29]{Fux}, which is in the Dorian mode. The integer $7$ is the \textit{interval} between the notes d and a, corresponding to $2$ and $9$ respectively. The linear polynomial $2+7x$ denotes the associated \textit{contrapuntal interval}. The transition from a bar to the next is a \textit{progression}.}
 \label{Fux}
\end{figure}

The simplest case of first-species counterpoint consists of two voices, \textit{cantus firmus} and \textit{discantus}, whose notes occur simultaneously and have the same duration. Thus, given a note in the cantus firmus, the corresponding one in the discantus is determined by the interval between them. Each interval is a consonance and the composition must satisfy certain rules. See Figure~\ref{Fux}.

The following is the codification of the \textit{strict style} given in \cite{MazzMuzz}, which is based on \cite{Tittel}. Here, for accessibility of the bibliographical sources, we justify each rule following the standard counterpoint books  \cite{Fux,Jep}.

The strict style consists of the following \textit{data}:
\begin{itemize}
\item \textit{Pitch space}. We will work with the usual group $\mathbb{Z}_{12}$ of pitch classes of the equal-tempered scale\footnote{For simplicity, we assume that it extends endlessly in both directions. The choice of this particular scale is only illustrative, but irrelevant. It could equally denote a scale in a different tuning.} $\mathbb{Z}$. We denote by $\pi$ the natural projection from $\mathbb{Z}$ to $\mathbb{Z}_{12}$.
\item \textit{Intervals}. The intervals between pitch classes and pitches also correspond to the groups $\mathbb{Z}_{12}$ and $\mathbb{Z}$, since they are just differences.
\item \textit{Diatonic scale}. The basic pattern of modes is the subset $X$ of $\mathbb{Z}_{12}$, where $X=\{0,2,4,5,7,9,11\}$, and its extension $\pi^{-1}(X)$ to $\mathbb{Z}$.
\item \textit{Contrapuntal intervals}. Each interval between a given note $c$ in the cantus firmus and its discantus $d$ corresponds uniquely to a linear polynomial $c+(d-c)x$. This expression is equivalent to giving the two voices explicitly as a pair $(c,d)$, but our intention is \textit{to stress the role of intervals}. We denote by $P_1$ the additive Abelian \textit{group of contrapuntal intervals}.

\item \textit{Consonances and dissonances}. We define consonances (unison, thirds, fifth, sixths) and dissonances (seconds, fourth, sevenths) by \[K=\{0,3,4,7,8,9\}\text{ and }D=\{1,2,5,6,10,11\},\] respectively. Thus, $\{K,D\}$ is a \textit{partition} of the intervals group $\mathbb{Z}_{12}$.

The consonance/dissonance partition of $\mathbb{Z}$ is the naturally induced by $\{K,D\}$, namely $\{\pi^{-1}(K),\pi^{-1}(D)\}$. In turn, by classifying contrapuntal intervals $c+(d-c)x$ according to whether $d-c$ is a consonance or a dissonance, we have a partition of the contrapuntal intervals group $P_1$ into \textit{contrapuntal consonances and dissonances}.
\end{itemize}

\textit{A composition of first-species counterpoint} is a finite sequence $\xi_1,\xi_2,\dots, \xi_n$ of contrapuntal intervals essentially satisfying the following \textit{rules}. 

\subsection{Preliminary rules}\label{prerul}
\begin{itemize}
\item Each contrapuntal interval is a contrapuntal consonance \cite[p.~27]{Fux}.
\item We work with consonances up to the tenth \cite[p.~	112, 5]{Jep}, that is, for each contrapuntal interval $c+kx$, $k\in\{0,3,4,7,8,9,12,15,16\}$.
\item Both $c$ (cantus firmus) and $d$ (discantus $c+k$) of each contrapuntal interval $c+kx$ are in the diatonic scale.\footnote{Moreover, the choice of a distinguished element $t$ in $X$ determines a \textit{mode} with \textit{tonic} $t$, where the composition occurs. Dorian, Phrygian, Mixolydian, Eolian, and Ionian modes are determined by the tonics $2$, $4$, $7$, $9$, and $0$, respectively. The tonic is essentially stressed in the first and last contrapuntal intervals of the composition \cite[p.~31]{Fux}, but the role of modes in progressions is secondary since the latter do not depend on a particular note (Section~\ref{strtrans}). See \cite[pp.~59-82]{Jep} for a musical introduction to modes.} 
\item Given a \textit{progression} $(c+kx,c'+k'x)$, from a contrapuntal consonance to its successor in a composition, the maximum change of a voice is an octave \cite[p.~109]{Jep}. Formally, $|c'-c|\leq 12$ and $|d'-d|\leq 12$, where $d=c+k$ and $d'=c'+k'$.
\end{itemize}

\subsection{Progression rules}\label{prorul}

According to \cite{MazzMuzz,Fux,Jep}, we divide progressions $(c+kx,c'+k'x)$ into three categories. 
\begin{itemize}
\item \textit{Inadmissible}: 
\begin{itemize}
\item[]  \textit{Unison repetitions} \cite[p.~	112, 3]{Jep}.
\item[] \textit{Parallel perfect consonances} \cite[p.~22]{Fux}. Those progressions with $k\in \{0,7,12\}$, $k=k'$, and $c\neq c'$.
\item[] \textit{Hidden parallel perfect consonances} \cite[p.~22]{Fux}. Those satisfying $k'\in\{0,7,12\}$, $(d'-d)(c'-c)>0$, and $k\neq k'$.
\item[] \textit{Tritones} \cite[p.~35]{Fux}. They satisfy $|c'-c|=6$ or $|d'-d|=6$.
\item[] \textit{Too large skips} \cite[p.~27, Footnote~1]{Fux}. If $7<|c'-c|<12$ or $7<|d'-d|<12$. The octave is regarded as a sort of repetition \cite[p.~112, 7]{Jep} and is accepted.
\end{itemize}
\item \textit{Bad}: Those that are not inadmissible but fall into some of the following cases. Bad progressions are not strictly avoided.
\begin{itemize}
\item[] \textit{Imperfect consonances by similar skips}\footnote{It is a weak form of the combination of the rules \textit{Parallel imperfect consonances} and \textit{Hidden parallel imperfect consonances}.} \cite[p.~112, 7]{Jep}. This means that $k'\in \{3,4,8,9,15,16\}$, $(d'-d)(c'-c)>0$, $|c'-c|>2$, $|d'-d|>2$, and $5<|c'-c|< 12$ or $5<|d'-d|<12$.
\item[] \textit{Hidden tritones} \cite[p.~35, Footnote~9]{Fux}. If $d'-c\equiv 6 \pmod{12}$ or $d-c'\equiv 6 \pmod{12}$.
\end{itemize}
\item \textit{Good}: All other progressions.
\end{itemize}

\subsection{Strict style modulo translation}\label{strtrans}

\textit{Since all rules can be expressed in terms of intervals}, they are invariant under translation of the progressions. Thus, to understand them, it is enough to study the representatives under translation of the form
\[(0+kx,c'+k'x).\]
See Table~\ref{tab:progtrans} for the counting of all progressions and their types. These outcomes were obtained with Code~1 in the Online Supplement.
\begin{table}
\begin{center}
\begin{tabular}{|c|c|c|}
\hline
\multirow{9}{*}{$1057$ prog.}&\multirow{6}{*}{$671$ inadmissible} & $22$ parallel fifths  \\
\cline{3-3}
 & & $49$ parallel eights and unisons \\
\cline{3-3}
& & $88$ hidden fifths  \\
\cline{3-3} 
 & & $128$ hidden eights and unisons  \\
\cline{3-3} 
 & & $170$ tritones  \\
\cline{3-3} 
& & $434$ too large skips  \\
\cline{2-3} & \multirow{2}{*}{$64$ bad} & $38$ imp. cons. by sim. skips  \\
\cline{3-3}& & $26$ hidden tritones \\
\cline{2-3} & \multicolumn{1}{c}{$322$ good}   &\\
\hline
\end{tabular}
\end{center}
\caption{Counting of the strict style progressions modulo translation and their types.}
\label{tab:progtrans}
\end{table}

\section{The principles of Mazzola's model}\label{modprin}

The following are the conceptual principles of the model; the mathematical details are in Section~\ref{mod}.

\subsection{Counterpoint rules obey mathematical laws based on symmetries}

Rather than an acoustic or psychological theory, counterpoint is a \textit{composition theory} with its own \textit{logic}. The fourth is an acoustic consonance, but it is a dissonance in counterpoint. Similarly, to justify the tritone prohibitions by saying that these intervals are hard to sing is not satisfactory \cite[p.~75]{Beau} given the versatility of instrumental counterpoint.

On the other hand, \textit{symmetries} (for example, transpositions and inversions) are \textit{the natural operations that occur in counterpoint}, so they are a reasonable basis for a model.

\subsection{We restrict to the intervals modulo octave}

The model does not intend to explain Fux counterpoint rules, but their reduction modulo octave, in the sense of Section~\ref{redstrsty}. This choice is just a \textit{simplification procedure}. 

\subsection{The division of intervals into consonances and dissonances has a mathematical conceptual characterization}

The partition of the twelve intervals modulo octave into consonances (unison, thirds, fifth, and sixths) and dissonances (seconds, fourth, sevenths) has a \textit{unique symmetry} that interchanges them. This fundamental property, discovered by Mazzola, \textit{characterizes} this partition together with a monoid property for consonances \cite[p.~76]{Beau}.

\subsection{The discantus is a tangential alteration of the cantus firmus}\label{tang}
First, there is an analogy of the intervals group $\mathbb{Z}_{12}$ with a \textit{differentiable manifold}, with associated tangent spaces, given the characterization of this Abelian group as the product $\mathbb{Z}_3\times \mathbb{Z}_4$, which can be interpreted as a \textit{torus of thirds} \cite[Section~12.1]{Beau}.

On the other hand, \textit{the discantus is a variation or alteration of the cantus firmus}, according to \cite[p.~73]{Beau} and \cite[p.~549]{Inicount}, and \textit{alterations are tangents} \cite[Section~7.5]{MazzolaTopos}. Thus, \textit{dual numbers}, which are the natural tangents in algebraic geometry \cite[p.~80]{Hartshorne}, suitably model this situation. 

This point of view of alterations as tangents is likely inspired by infinitesimal deformations in algebraic geometry \cite[Example~9.13.1]{Hartshorne}, which are related to dual numbers. 

\subsection{Alternation and deformation}

Contrapuntal tension is not only vertical (between cantus firmus and discantus) but horizontal (between an interval and its successor), as reflected by the division of consonances into perfect (unison/octave and fifth) and imperfect ones (thirds and sixths). This suggest a trace of dissonance inside consonance. These ideas are inspired by the work of Klaus-Jürgen Sachs as explained in \cite[Section~14.1]{Beau} and \cite[p.~646]{MazzolaTopos}.

To translate the idea of horizontal tension in the model, we consider progressions of consonances (interval to successor) that can be regarded as symmetry-deformed progressions from a dissonance to a consonance \cite[Section~31.1]{MazzolaTopos}. This \textit{alternation} (or contrast) process resembles the musical resolution sense of dissonances into consonances.

\section{Summary of mathematical counterpoint theory}\label{mod}

The following is a synthesis of a more robust counterpoint theory \cite{Theo}. We start with the following basic data.

\begin{itemize}
\item The \textit{ring} $\mathbb{Z}_{12}$, which corresponds to the twelve intervals (up to the octave).

\item The group of \textit{symmetries} of $\mathbb{Z}_{12}$, denoted by $\sym(\mathbb{Z}_{12})$, which consists of all affine automorphisms of the form $e^ab:\mathbb{Z}_{12}\longrightarrow \mathbb{Z}_{12}:r\mapsto br+a$, where $b\in \mathbb{Z}_{12}^*=\{1,3,5,7,11\}$. The symmetries of the form $e^a1$, or $e^a$ for short, called \textit{translations}, correspond to transposition in music.
\item The \textit{partition} $\{K,D\}$ of $\mathbb{Z}_{12}$ into \textit{consonances} $K$ and \textit{dissonances} $D$, where $K=\{0,3,4,7,8,9\}$ and $D=\{1,2,5,6,10,11\}$. Mazzola's fundamental observation is that $e^2 5$ is the \textit{unique} $p\in \sym(\mathbb{Z}_{12})$ such that $p(K)=D$. 
\item The \textit{dual numbers ring} $\mathbb{Z}_{12}[\e]$, which models contrapuntal intervals modulo octave. It consists of all linear polynomials $a+b\e$ subject to the relation $\e^2=0$. Formally, it is the quotient ring $\mathbb{Z}_{12}[x]/\left\langle x^2\right\rangle$, where $\e$ is the class of $x$.
\item The group of \textit{symmetries} of $\mathbb{Z}_{12}[\e]$, denoted by $\sym(\mathbb{Z}_{12}[\e])$, which consists of all affine automorphisms of the form $e^{u+v\e}(c+d\e):\mathbb{Z}_{12}[\e]\longrightarrow \mathbb{Z}_{12}[\e]:x+y\e\mapsto (c+d\e)(x+y\e)+(u+v\e)$, where $c\in \mathbb{Z}_{12}^*$. 
\item The induced partition $\{K[\e],D[\e]\}$ of $\mathbb{Z}_{12}[\e]$ into \textit{contrapuntal consonances} $K[\e]$ and \textit{contrapuntal dissonances} $D[\e]$, where 
\[Y[\e]=\{r+k\e\ | \ r\in \mathbb{Z}_{12} \text{ and }k\in Y \}\]
for $Y=K,D$.
\end{itemize}

We consider pairs $(\xi,\eta)$ of contrapuntal consonances in $K[\e]$, which we call \textbf{progressions}, and aim to determine when they are valid for counterpoint. The three principles in Sections~\ref{alternation}, \ref{locchar}, and \ref{variety} serve this purpose.

\subsection{Alternation and repetition}\label{alternation}

The model deals with \textbf{polarized progressions}, that is, progressions $(\xi,\eta)$ such that $\xi\in g(D[\e])$ and $\eta\in g(K[\e])$ for some $g\in\sym(\mathbb{Z}_{12}[\e])$. As proved in \cite{Theo}, \textit{polarized progressions are all but repetitions}. Since the model selects among them those that are optimal in a musical sense (Definition~\ref{def}), we say
that \textit{the model does not decide on the nature of repetitions}.

\subsection{Local characterization of consonances and dissonances}\label{locchar}

Although alternation helps to regard $\xi$ as a deformed dissonance and $\eta$ as a deformed consonance, we should ensure that they resemble dissonances and consonances in an authentic way. Thus, we would want a property for the partition $\{K[\e],D[\e]\}$, analogous to the uniqueness one of $\{K,D\}$, that defines contrapuntal consonances and dissonances.

The symmetry $e^{2\e} 5$ of $\mathbb{Z}_{12}[\e]$ is a natural extension of $e^2 5$ sending $K[\e]$ to $D[\e]$ since $e^{2\e}5$ is simple and acts on the interval part $d$ of a dual number $a+b\e$ just as $e^25$. But in this case \textit{it is not the unique} sending $K[\e]$ to $D[\e]$. For example, $e^{2\e+1} 5$ also does. 

However, we have the following \textit{local uniqueness property} \cite[Proposition~51]{MazzolaTopos}, which is a sort of \textit{characteristic one} \cite[Section~14.2]{Beau} of the consonance/dissonance partition $\{z+K\e,z+D\e\}$ of the \textit{fiber} $z+\mathbb{Z}_{12}\e$, thus offering a local definition of contrapuntal consonance and dissonance.\footnote{
See \cite[Section~4.2]{Theo} for the emergence of a similar \textit{characterizing} property of the consonance/dissonance partition $\{z+K\e,z+D\e\}$. The two properties are related in \cite[Section~10]{Theo}, where there is a proof that they are equivalent for the deformed partition $\{g(K[\e]),g(D[\e])\}$ at a given fiber \cite[Theorem~10.6]{Theo}.} 

\begin{theorem}\label{indxdich}
For each cantus firmus note $z\in \mathbb{Z}_{12}$, the symmetry $p^z[\e]$, defined by $p^z[\e]=e^z\circ e^{2\e}5\circ e^{-z}=e^{8z+ 2\e}5$, is the unique in $\sym(\mathbb{Z}_{12}[\e])$ that
\begin{itemize}
\item[1.] leaves invariant the fiber $z+\mathbb{Z}_{12}\e$ and
\item[2.] sends $K[\e]$ to $D[\e]$.
\end{itemize}
In particular, 
\begin{equation*}
p^z[\e](z+K\e)=z+D\e.
\end{equation*}
\end{theorem}

If we apply this property to the deformed partition $\{g(K[\e]),g(D[\e])\}$, we obtain the following one.
\begin{itemize}
\item The symmetry $p^z[\e]$ is the unique $p'\in \sym(R[\e])$ that leaves invariant the fiber $z+\mathbb{Z}_{12}\e$ and sends $g(K[x])$ to $g(D[x])$.
\end{itemize}
By \cite[Proposition~10.5]{Theo}, this condition is equivalent to the following equation.
\begin{equation}\label{two}
p^z[\e](g(K[\e]))=g(D[\e])
\end{equation}

Mazzola requires it for the cantus firmus $z$ of $\xi$, which ensures that $\xi$ is a local dissonance.

\subsection{Variety}\label{variety}

In this model, the variety principle of counterpoint \cite[p.~21]{Fux} corresponds to the condition that there is a maximum of alternations from $\xi$, that is, the cardinality of $g(K[\e])\cap K[\e]$ is maximum among all $g\in \sym(R[\e])$ such that 1. $\xi \in g(D[\e])\cap K[\e]$ (alternation) and 2. Equation~\eqref{two} holds for the cantus firmus $z$ of $\xi$ (local dissonance).

\subsection{Admitted successors}

Now we can give the definition of admitted successor of a contrapuntal interval. 

\begin{figure}
\centering
\def\svgscale{0.8}
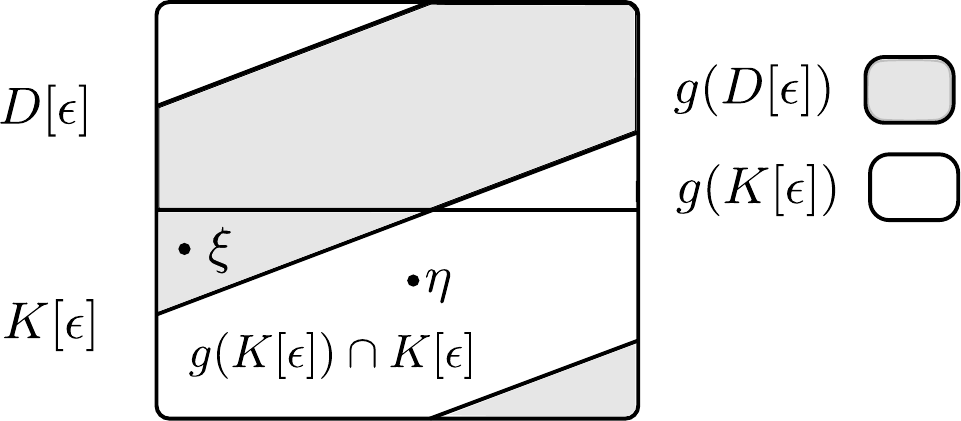
\caption{Here, $g$ is a deformation symmetry and $\eta$ an admitted successor of $\xi$.}
 \label{counterpoint}
\end{figure}

\begin{definition}\label{def} 
A \textbf{contrapuntal symmetry} for a consonance $\xi\in K[\e]$, where $\xi=z+k\e$, is a symmetry $g$ of $\mathbb{Z}_{12}[\e]$ such that 
\begin{enumerate}
\item $\xi\in g(D[\e])$,
\item the symmetry $p^z[\e]$ sends $g(K[\e])$ to $g(D[\e])$, and 
\item the cardinality of $g(K[\e])\cap K[\e]$ is maximum among all $g$ satisfying 1 and 2. 
\end{enumerate}
Note that the contrapuntal symmetry for a given consonance is not required to be unique. 

An \textbf{admitted successor} of a consonance $\xi\in K[\e]$ is an element $\eta$ of $g(K[\e])\cap K[\e]$ for some contrapuntal symmetry $g$. See Figure~\ref{counterpoint}. 

If $\eta$ is an admitted successor of $\xi$, we say that the progression $(\xi,\eta)$ is \textbf{allowed}. If it does not happen and $(\xi,\eta)$ is polarized, we say that it is \textbf{forbidden}.
\fini
\end{definition}

\subsection{Admitted successors computation}\label{adsuccomp}

Denote by $H$ the group\footnote{In fact, it is a group since it consists of all symmetries that leave invariant the fiber $K\e$ at $0$.} of all symmetries of $\mathbb{Z}_{12}[\e]$ of the form $e^{v\e}(c+d\e)$.

According to \cite{Theo}, after several simplifications, the admitted successors of $z+k\e\in K[\e]$ are the elements of the sets of the form
\begin{equation}\label{transpose}
e^z(h(K[\e])\cap K[\e]),
\end{equation}
where $h=e^{v\e}(c+d\e)\in H$ and 
\begin{enumerate}
\item $v\in k-cD$, 
\item $5v+2=c2+v$, and 
\item the value $\rho\sum\limits_{i=0}^{\rho-1}|K_i||K_{e^vc(i)}|$ (the cardinality of $h(K[\e])\cap K[\e]$) is maximum among all $h=e^{v\e}(c+d\e)\in H$ satisfying 1 and 2, where $K_{i}=\{k\in K\ |\ k \equiv i\pmod{\rho}\}$, $\rho=\gcd (d,12)$, and $e^vc$ is reduced modulo $\rho$. 
\end{enumerate}
The symmetries $h$ satisfying the previous conditions are also contrapuntal symmetries for $k\e$.

This means that, to compute the admitted successors of a consonance $z+k\e$, it is enough to do so for $k\e$, and then apply the translation $e^z$ (Equation~\eqref{transpose}). This agrees with the previous result that counterpoint rules do not depend on particular notes; see Section~\ref{strtrans}. For these reasons, from now on \textit{we assume that all progressions are of the form $(k\epsilon, c'+k'\epsilon)$}. 

Finally, the counterpoint symmetries and admitted successors of $k\e\in K[\e]$ can be computed with \textit{Hichert's algorithm}. This algorithm ranges over all symmetries $h$ satisfying (1) and (2), and updates in each step the set of those whose associated cardinalities, according to (3), are maximum so far. Then, with the counterpoint symmetries at hand, we compute the associated successors sets according to the formula
\begin{equation}\label{presuc}
e^{v\e}(c+d\e)(K[\e])\cap K[\e]=\bigsqcup\limits_{r\in \mathbb{Z}_{12}}cr+((cK+v+dr)\cap K)\e.
\end{equation}
 
The Table~1 in the Online Supplement shows all counterpoint symmetries and successors. It was obtained with the Code~2 in the Online Supplement.  

According to Code~3 in the Online Supplement, out of the $287$ progressions that occur in the diatonic scale $X$, $6$ of them are repetitions (non-polarized), $250$ are allowed, and the remaining $31$ are forbidden. Next, we reduce the strict style so as to compare it with the model.  

\section{The reduced strict style}\label{redstrsty}

Here we reduce, modulo $12$ (up to the octave), the strict style modulo translation studied in Section~\ref{strsty}. We use the projection that sends a progression $(kx,c'+k'x)$ in the strict style to $(\pi(k)\e,\pi(c')+\pi(k')\e)$, where $\pi:\mathbb{Z}\longrightarrow \mathbb{Z}_{12}$ is the natural projection. 

We first note that \textit{this projection covers all progressions in a diatonic scale}. In fact, each such progression $(k\e,c'+k'\e)$, regarded as $(kx,c'+k'x)$, satisfies all preliminary rules in Section~\ref{prerul}, except perhaps the condition that the maximum change between the discantus notes is an octave. In such a case, we transpose the second interval an eight downwards and obtain an strict style progression that projects on $(k\e,c'+k'\e)$.

Now we define the progression \textit{rules of the reduced strict style}. According to \cite{MazzMuzz}, a progression in a diatonic scale is
\begin{itemize}
\item[] \textbf{good}, if it is the projection of at least one good progression,
\item[] \textbf{inadmissible}, if it is derived from nothing but inadmissible progressions, and
\item[] \textbf{bad}, if it is the projection of at least one bad progression, but not of a good one.
\end{itemize}

This definition leads to the following characterization of these rules. Code~4 in the Online Supplement contains some computations involved.

\subsection*{Inadmissible progressions}
Under projection, unison repetitions remain unchanged and they are good, parallel unisons or eights become parallel unisons or unison repetitions, parallel fifths become parallel fifths or fifth repetitions (good), and tritones become tritones possibly greater than the eight.

\begin{itemize}
\item[] \textbf{Parallel unisons and fifths} 

\textit{Some parallel unisons remain inadmissible}. If the skip is the tritone, then it is inadmissible. Now, parallel unisons come from other three kinds of progressions: 1. unison to octave, 2. parallel eights (inadmissible), and 3. octave to unison. The cases 1 or 3 are inadmissible if and only if $c'$ or $12-c'$ are greater than $7$, that is, $c'\in \{1,2,3,4,8,9,10,11\}$. Otherwise, we have good progressions by contrary motion and skips not too large. 

\textit{Parallel fifths are inadmissible}. They can only be the octave reduction of parallel fifths, which are inadmissible, because twelfths are not in the range of the strict style. 

\item[] \textbf{Tritones}

\textit{Tritones are inadmissible}. Tritones only come from tritones, which are inadmissible.

\item[] \textbf{Projected hidden parallelisms}

\textit{Hidden parallel fifths from a sixth are inadmissible}. They are only the projection of progressions of the same kind.

The remaining cases that have not tritones are not inadmissible (Code~4).

\item[] \textbf{Projected too large skips}

All cases fall into the previous inadmissible ones or are not inadmissible (Code~4).
\end{itemize}

\subsection*{Bad progressions}
Bad progressions in $\mathbb{Z}_{12}[\e]$ have to be projections of at least a bad progression in $\mathbb{Z}[\e]$.

\begin{itemize}

\item[] \textbf{Projected imperfect consonances by similar skips} 

\textit{Here, the unique bad progression is $(0+7\e,5+9\e)$}; see Code~4. The other ones are good.

\item[] \textbf{Hidden tritones}  

\textit{Hidden tritones are bad}. They only come from hidden tritones, which are bad.
\end{itemize}

To sum up, in the reduced strict style, the inadmissible progressions are all tritones and, excepting them, some parallel unisons, all parallel fifths, and all hidden parallel fifths from a sixth. Moreover, as noted in \cite{MazzMuzz}, \textit{only the parallel fifths and tritone rules preserve their generality} (unrestricted validity).

\subsection*{On semantics}

However, out of the previous categories (\textit{inadmissible}, \textit{bad}, \textit{good}), the most appropriate, to provide a quantitative assessment of the agreement of the model with the original rules, seem to be the \textit{inadmissible} and \textit{bad} ones. Certainly, an inadmissible progression only comes from inadmissible progressions so it is \textit{unequivocally inadmissible} and bad progressions are optional, and unequivocally bad in this case (check). In contrast, good progressions also come from inadmissible or bad progressions, except for the four imperfect consonance repetitions, which only come from good ones, as we deduce from the following definition.  

We have a \textbf{refined semantics} of the reduced strict style. We preserve the definition of inadmissible and bad progressions. Regarding a remaining good progression, we sub-classify it as
\begin{itemize}
\item[] \textbf{good-good}, if it is derived from nothing but good progressions,
\item[] \textbf{ambiguous}, if it also comes from at least an inadmissible progression, and
\item[] \textbf{good-bad}, if it also comes from at least a bad progression but not an inadmissible one.
\end{itemize}
This definition offers a fairer classification: ambiguous progressions have not a defined value comparable to \textit{allowed} and \textit{forbidden}.

Table~\ref{tab:progredu} contains the counting of all reduced strict style progressions and their types, according to Code~4 in the Online Supplement. 
\begin{table}
\begin{center}
\begin{tabular}{|c|c|c|}
\hline
\multirow{9}{*}{$287$ prog.}&\multirow{4}{*}{$74$ inadmissible} & $10$ parallel fifths*  \\
\cline{3-3}
 & & $9$ parallel unisons \\
\cline{3-3}
& & $13$ hidden fifths from a sixth \\
\cline{3-3} 
 & & $45$ tritones*  \\
\cline{2-3} & \multirow{2}{*}{$23$ bad} & $1$ proj. imp. cons. by  sim. skips  \\
\cline{3-3}& & $22$ hidden tritones \\
\cline{2-3} & \multirow{3}{*}{$190$ good}   & $4$ good-good  \\
\cline{3-3}& & $16$ good-bad \\
\cline{3-3}& & $170$ ambiguous \\
\hline
\end{tabular}
\end{center}
\caption{Counting of the reduced strict style progressions and their types. The symbol * refers to rules whose generality remains under projection.}
\label{tab:progredu}
\end{table}

\section{The reduced strict style and the model}\label{comp}

Table \ref{tab:comp} compares all types of strict style progressions with allowed and forbidden ones, according to Code~5 in the Online Supplement. 

\begin{table}
\begin{center}
\begin{tabular}{|c|c|c|c|c|c|c|c|c|}
\hline
  & inad.& inad.*  & bad & good & good* & good-good & good-bad & amb.\\
\hline
allowed & $55$ & $36$& $17$ & $178$ & $197$& $0$ & $16$&$162$\\
\hline
forbidden & $19$ & $19$ & $6$ & $6$ & $6$& $0$& $0$& $6$\\
\hline
repetitions & $0$ & $0$ & $0$ & $6$  &  $6$ & $4$& $0$& $2$\\
\hline
\end{tabular}
\end{center}
\caption{Inadmissible, bad, and good progressions (reduced strict style) versus allowed and forbidden ones (model). The symbol * refers to rules whose generality remains under projection.}
\label{tab:comp}
\end{table}

Table \ref{tab:inadbad} shows all allowed and forbidden kinds of inadmissible and bad progressions, according to Code~6 in the Online Supplement. We observe that \textit{the model predicts the parallel fifths prohibition} and some tritone rules.
\begin{table}
\begin{center}
\begin{tabular}{|c|c|c|c|c|c|c|}
\hline
  &par. 5ths& par. un. & hid. 5ths& trit.& $(0+7\e,5+9\e)$& hid. trit.\\
\hline
allow. & $0$& $8$ & $13$& $36$ & yes& $16$\\
\hline
forb. &  $10$ & $1$ & $0$& $9$ &  no & $6$ \\
\hline
\end{tabular}
\end{center}
\caption{Allowed and forbidden kinds of inadmissible and bad progressions.}
\label{tab:inadbad}
\end{table}

We can measure the agreement between the model and the reduced style by means of the number of matches and mismatches. The matches are all allowed good, forbidden inadmissible, and (possibly) forbidden bad progressions. The mismatches are all forbidden good and allowed inadmissible ones. From Table~\ref{tab:comp}, we observe that there are $203$ matches, including forbidden bad progressions, and $61$ mismatches. If we only take into account rules whose generality remains under projection, we obtain $222$ matches and $42$ mismatches. 

We define matches and mismatches in the refined semantics by replacing good progressions for good-good progressions in the previous paragraph. In the case of the previous model, there are $25$ matches and $55$ mismatches. But, since there are only six good-good progressions, these measures are not appropriate. For example, the trivial model that forbids all progressions, except the four imperfect consonance repetitions, has $101$ matches and no mismatches---the best possible result. This means that the ambivalence of $186$ good progressions (ambiguous and good-bad) does not allow an appropriate semantics on reduced progressions that captures the essence of the original good ones.\footnote{In other words, a strong attribute of being good is invisible or undetectable from the point of view of the reduced progressions.} However, we use these measures to discard possible ambiguities.

Thus, the careful review of the semantics offers a different point of view on the results of the model. For instance, we observe that the twelve examples (other than repetitions) of progressions forbidden by Mazzola but not by Fux, according to \cite[Fig.~3b]{Tympolem}, are the six forbidden bad together with the six forbidden good progressions in Table~\ref{tab:comp}. The last six progressions are ambiguous according to the refined semantics, so the twelve examples are not actually wrong predictions of the model. 

\section{Some questions on the model}

The conceptual review of the principles in Sections~\ref{tang} and \ref{locchar}   leads to the following questions.

\begin{itemize}
\item Why is the discantus necessarily a tangential alteration of the cantus firmus?
\item Why do we not require the local characterization condition on $\{g(K[\e]),g(D[\e])\}$ for the cantus firmus of a possible successor $\eta$?
\item More radically, why do we not require the local characterization on $\{g(K[\e]),g(D[\e])\}$ for all fibers?
\end{itemize}
These questions inspire the following three variations of the model.

\section{First variation of the model}\label{natvar}

Besides the uniqueness condition on the partition $\{K,D\}$, the dual numbers structure on counterpoint intervals is the most important structural feature of Mazzola's model. However, it is not utterly clear why we use the structure of the dual numbers ring to model intervals and not another one. 

Another simple choice for a ring structure on counterpoint intervals is that induced by the \textit{product} ring $\mathbb{Z}_{12}\times \mathbb{Z}_{12}$, whose elements $(c,d)$ can be regarded as a pair of cantus firmus and discantus notes that occur simultaneously. If we transfer its structure to contrapuntal intervals, under the map that sends $(c,d)$ to $(c,d-c)$ (cantus firmus and interval with the discantus), we obtain (check) the structure of the quotient ring $\mathbb{Z}_{12}[x]/\left\langle x^2-x\right\rangle$---a variation of the dual numbers ring. In this case, we denote the class of $x$ by $\x$, and the ring by $\mathbb{Z}_{12}[\x]$. The latter consists of linear polynomials $a+b\x$, where $a,b\in\mathbb{Z}_{12}$ and $\x^2=\x$, that is, $\x$ is idempotent.

The rings $\mathbb{Z}_{12}[\e]$ and $\mathbb{Z}_{12}[\x]$ are equal as Abelian groups but they differ in their products. In $\mathbb{Z}_{12}[\e]$, $d\e d'\e=dd'\e^2=0$ for all $d,d'\in \mathbb{Z}_{12}$, so intervallic variations are regarded as \textit{infinitesimals}, whereas in $\mathbb{Z}_{12}[\x]$, $d\x d'\x=dd'\x ^2=dd'\x$, so these variations are just \textit{integers modulo $12$}.

According to \cite{Theo}, like in the case of dual numbers, symmetries of $\mathbb{Z}_{12}[\x]$ are of the form $e^{u+v\x}(c+d\x)$, but both $c$ and $c+d$ are in $\mathbb{Z}_{12}^*$ now. In the same way, the extension of consonances and dissonances $\{K[\x],D[\x]\}$ to counterpoint intervals can be defined and Theorem~\ref{indxdich} remains valid by changing $\e$ for $\x$. The non-polarized symmetries, according to the definition in Section~\ref{alternation}, are all repetitions together with all parallelisms by tritone leaps in this case. Definition~\ref{def} is the same because the motivations remain valid, except that we replace $\e$ by $\x$. Regarding Section~\ref{adsuccomp}, now $H$ consists of all symmetries of $\mathbb{Z}_{12}[\x]$ of the form $e^{v\x}(c+d\x)$, which satisfy $c,c+d\in \mathbb{Z}_{12}^*$. However, conditions 1 and 2 become 
\begin{enumerate}
\item $v\in k-(c+d)D$ and
\item $5v+2=(c+d)2+v$, 
\end{enumerate}
whereas 3 remains unchanged. We compute the successors sets with the formula:
\begin{equation}
e^{v\x}(c+d\x)(K[\x])\cap K[\x]=\bigsqcup\limits_{r\in \mathbb{Z}_{12}}cr+(((c+d)K+v+dr)\cap K)\x.
\end{equation}

The results of this model, and their comparison to the strict reduced style, are summarized in Tables \ref{tab:comp1} and \ref{tab:inadbad1}. Here, out of the $287$ progressions that occur in a diatonic scale, $7$ are non-polarized, $240$ are allowed, and $40$ are forbidden. These results correspond to Code~7 in the Online Supplement. We observe that the model predicts $21$ inadmissible progressions. Actually, it predicts $18$ out of the $19$ inadmissible progressions of the classical model together with three new hidden fifths from a sixth. We lose the parallel unison by tritone skip, since it is non-polarized here. Regarding bad progressions, it predicts the same six hidden tritones. However, in this case we have $198$ matches and $65$ mismatches with the reduced strict style. In the refined semantics the respective measures are $27$ and $52$.


\begin{table}
\begin{center}
\begin{tabular}{|c|c|c|c|c|c|c|c|c|}
\hline
  & inad.&  bad & good  & good-good & good-bad & amb.\\
\hline
allowed & $52$ &  $17$ & $171$ & $0$ & $15$&$156$\\
\hline
forbidden & $21$ & $6$ & $13$ &  $0$& $1$& $12$\\
\hline
non-polarized & $1$ & $0$ & $6$  &  $4$& $0$& $2$\\
\hline
\end{tabular}
\end{center}
\caption{Inadmissible, bad, and good progressions (reduced strict style) versus allowed and forbidden ones (first variation). }
\label{tab:comp1}
\end{table}

\begin{table}
\begin{center}
\begin{tabular}{|c|c|c|c|c|c|c|}
\hline
  &par. 5ths& par. un. & hid. 5ths & trit.& $(0+7\x,5+9\x)$& hid. trit.\\
\hline
allow. & $0$& $8$ & $10$& $36$ & yes& $16$\\
\hline
forb. &  $10$ & $0$ & $3$& $8$ &  no & $6$ \\
\hline
non-pol. &  $0$ & $1$ & $0$& $1$ &  no & $0$ \\
\hline
\end{tabular}
\end{center}
\caption{Allowed and forbidden kinds of inadmissible progressions (first variation).}
\label{tab:inadbad1}
\end{table}

\section{Second variation of the model}\label{2var}

If we also require the local characterization on the fibers of the possible successors of a consonance, then we change condition (3) in Definition~\ref{def} for the following one, where $P(z)$ denotes the property (2) in Definition~\ref{def}.
\begin{itemize}
\item The cardinality of $\{z'+k'\e \in g(K[\e])\cap K[\e]\ |\ P(z')\}$ is maximum among all $g$ satisfying (1) and (2). 
\end{itemize}

Remarkably, this variation is equivalent to the following one, as proved in \cite[Section~11.3]{Theo}.

\section{Third variation of the model}\label{3var}

If we require the local characterization property on the deformed partition for all fibers, which amounts to a \textit{global condition}, then we change condition (2) in Definition~\ref{def} for the following one.
\begin{itemize}
\item For all $z\in \mathbb{Z}_{12}$, $P(z)$ holds.
\end{itemize}

According to \cite[Section~11.3]{Theo}, so as to compute the admitted successors according to this definition, it is enough to follow the procedure in Section~\ref{adsuccomp} by adding to condition (2) the equation $5d=d$. The final results of the second and third variations coincide and they are the same for the case of $\mathbb{Z}_{12}[\x]$ \cite[Sections~11.2-11.3]{Theo}, except for the difference that the inadmissible parallel unison by tritone skip is forbidden in the dual numbers case but non-polarized in the case of $\mathbb{Z}_{12}[\x]$.

The comparison of these results with the reduced strict style are in Tables~\ref{tab:comp2} and \ref{tab:inadbad2}. They correspond to Code~8 in the Online Supplement. We only put the results for the case of $\mathbb{Z}_{12}[\e]$. Here, there are six repetitions, $46$ forbidden progressions, and $235$ allowed ones, that occur in a diatonic scale.

\begin{table}
\begin{center}
\begin{tabular}{|c|c|c|c|c|c|c|c|c|}
\hline
  & inad.&  bad & good  & good-good & good-bad & amb.\\
\hline
allowed & $53$ &  $15$ & $167$ & $0$ & $12$&$155$\\
\hline
forbidden & $21$ & $8$ & $17$ &  $0$& $4$& $13$\\
\hline
repetitions & $0$ & $0$ & $6$  &  $4$& $0$& $2$\\
\hline
\end{tabular}
\end{center}
\caption{Inadmissible, bad, and good progressions (reduced strict style) versus allowed and forbidden ones (final variations).}
\label{tab:comp2}
\end{table}

\begin{table}
\begin{center}
\begin{tabular}{|c|c|c|c|c|c|c|}
\hline
  &par. 5ths& par. un. & hid. 5ths& trit.& $(0+7\e,5+9\e)$& hid. trit.\\
\hline
allow. & $0$& $4$ & $13$& $38$ & yes& $14$\\
\hline
forb. &  $10$ & $5$ & $0$& $7$ &  no & $8$ \\
\hline
\end{tabular}
\end{center}
\caption{Allowed and forbidden kinds of inadmissible and bad progressions (final variations).}
\label{tab:inadbad2}
\end{table}

The most important feature of these results, regarding the original model, is the \textit{prediction of new parallel unisons}. However, we have $196$ matches and $70$ mismatches. In the refined semantics, the respective measures are $29$ and $53$.

\section{Conclusions}\label{conc}

\subsection*{Economy}

Essentially, we deduce all mathematical results of the theory from just a fact: the uniqueness property of the consonance/dissonance partition. In particular, we deduce the parallel fifths prohibition. 

\subsection*{Generality and universality}
The models have a generalization to \textit{any} (not necessarily commutative) ring \cite{Theo}, taking the place of $\mathbb{Z}_{12}$. This is just an expression of the original intention of establishing an universal counterpoint by detecting the essential features of the Renaissance incarnation. Thus, the model paves the way to many other forms of counterpoint, which opens up a lot of fields of musical experimentation.

\subsection*{Quantitative interpretation of the results}
Table~\ref{tab:match} shows the matches and mismatches of all models with the reduced strict style. There is a number of mismatches in all cases, which might be interpreted as a weakness of the model, but the number of matches is about three times it.

The number of matches and mismatches can increase or decrease from the classical to the final variations according to the semantics used, but the differences do not exceed five units. 
\begin{table}
\begin{center}
\begin{tabular}{|c|c|c|c|c|c|c|c|c|}
\hline
  & matches &  mismatches & matches (ref. sem.)  & mismatches (ref. sem.)\\
\hline
classical & $203$ &  $61$ & $25$ & $55$ \\
\hline
first var. & $198$ & $65$ & $27$ &  $52$\\
\hline
final var. & $196$ & $70$ & $29$  &  $53$\\
\hline
\end{tabular}
\end{center}
\caption{Matches and mismatches of all models with the reduced strict style according to the original semantics and the refined semantics. Matches include all forbidden inadmissible, bad inadmissible, and allowed good progressions. Mismatches include all forbidden good and allowed inadmissible progressions. In the refined semantics, we change good progressions for good-good progressions, which do not come from inadmissible or bad progressions in the strict style.}
\label{tab:match}
\end{table}

\subsection*{Qualitative interpretation of the results}

Rather than a mere description of the Renaissance counterpoint rules (Section~\ref{strtrans}), we believe that the point is what mathematical counterpoint theory can explain about them. Mazzola's model explains the parallel fifths rule in its generality and some tritone prohibitions. Additionally, the first variation of the model explains some hidden fifths, and the final variations explain some parallel eights.

\subsection{Recovery of the original phenomenon from the model}\label{recov}

We can characterize perfect consonances as $0$ (which belongs to any ring) and the consonance with a general parallelism prohibition. The tritone can be characterized as the distinguished skip that makes a parallel unison forbidden. In the first variation, the tritone is the non-polarized skip.

In this way, we can recover the original rules by establishing as prohibitions the parallel and hidden perfect consonances, the distinguished skip, and the too large skips. 

\section{Further developments}\label{further}

\subsection*{Generalizations and structural approaches}

Further generalizations to categories, beyond that to arbitrary rings, could be developed. The use of the product ring $\mathbb{Z}_{12}\times \mathbb{Z}_{12}$ to model contrapuntal intervals, as an alternative to the dual numbers, helps express the model in terms of universal constructions in the category of rings, so it could be translated to other categories. 

Also, the mathematical counterpoint theory is an independent field of study. There are several open questions like whether there is an entirely structural proof of the \textit{Counterpoint theorem}, which establishes all counterpoint symmetries and admitted successors of a given consonance, beyond algorithmic computations. Some recent advances in that direction, like the counting formulas for successors sets cardinalities and a maximization criterion, can be found in \cite{Theo}.

\subsection*{Study of new counterpoint worlds}

The musical study of these generalizations is an open field of research. It can start by establishing rules based on the model, analogous to the original rules of counterpoint, following the procedure suggested in Section~\ref{recov}. Then, composition processes should generate new music in these worlds, as mentioned in Section~\ref{intro}. The counterpoint world induced by Scriabin's mystic chord has also been subject of study \cite{Scriabinworld}.

\subsection*{Musicology}

An initial comparison of the model's results with the \textit{Missa Papae Marcelli} can be found in \cite{Nieto}. The analysis of other works of Renaissance polyphony is also a pending task.   

\subsection*{Remaining species}

Although there are some advances regarding the second-species \cite{secsp}, a mathematical theory of the remaining species, and the cases of three or more voices, are not at hand.



\section*{Acknowledgements}
\addcontentsline{toc}{section}{Acknowledgements}

We thank Guerino Mazzola for his agreement to use the material \cite{MazzMuzz}. We also thank the two anonymous reviewers and Editor [name] for \dots

\section*{Funding}
\addcontentsline{toc}{section}{Funding}

This work was supported by Programa de Becas Posdoctorales en la UNAM 2019, which is coordinated by Direcci\'{o}n de Asuntos del Personal Acad\'{e}mico (DGAPA) at Universidad Nacional Aut\'{o}noma de M\'{e}xico. 

\section*{Supplemental online material}
\addcontentsline{toc}{section}{Supplemental online material}

Supplemental online material for this article can be accessed at \url{doi-provided-by-publisher} and/or \url{https://www.dropbox.com/s/1thilv2ji6cem0u/supplement.pdf?dl=0}. In the Online Supplement we include all \textit{Python} codes used in this paper.


\section*{Disclosure statement}
\addcontentsline{toc}{section}{Disclosure statement}

No potential conflict of interest was reported by the authors. \\ \\



\bibliographystyle{amsplain}
\bibliography{MySubmissionBibTexDatabase}

\addcontentsline{toc}{section}{References}

\end{document}

%% file: Fux.pdf_tex
\begingroup%
  \makeatletter%
  \providecommand\color[2][]{%
    \errmessage{(Inkscape) Color is used for the text in Inkscape, but the package 'color.sty' is not loaded}%
    \renewcommand\color[2][]{}%
  }%
  \providecommand\transparent[1]{%
    \errmessage{(Inkscape) Transparency is used (non-zero) for the text in Inkscape, but the package 'transparent.sty' is not loaded}%
    \renewcommand\transparent[1]{}%
  }%
  \providecommand\rotatebox[2]{#2}%
  \newcommand*\fsize{\dimexpr\f@size pt\relax}%
  \newcommand*\lineheight[1]{\fontsize{\fsize}{#1\fsize}\selectfont}%
  \ifx\svgwidth\undefined%
    \setlength{\unitlength}{436.93575516bp}%
    \ifx\svgscale\undefined%
      \relax%
    \else%
      \setlength{\unitlength}{\unitlength * \real{\svgscale}}%
    \fi%
  \else%
    \setlength{\unitlength}{\svgwidth}%
  \fi%
  \global\let\svgwidth\undefined%
  \global\let\svgscale\undefined%
  \makeatother%
  \begin{picture}(1,0.22477513)%
    \lineheight{1}%
    \setlength\tabcolsep{0pt}%
    \put(0,0){\includegraphics[width=\unitlength,page=1]{Fux.pdf}}%
    \put(-0.00074708,0.05381261){\color[rgb]{0,0,0}\makebox(0,0)[lt]{\lineheight{1.25}\smash{\begin{tabular}[t]{l}cantus firmus\end{tabular}}}}%
    \put(0.0239705,0.15641778){\color[rgb]{0,0,0}\makebox(0,0)[lt]{\lineheight{1.25}\smash{\begin{tabular}[t]{l}discantus\end{tabular}}}}%
    \put(0,0){\includegraphics[width=\unitlength,page=2]{Fux.pdf}}%
  \end{picture}%
\endgroup%

%% file: counterpoint.pdf_tex
\begingroup%
  \makeatletter%
  \providecommand\color[2][]{%
    \errmessage{(Inkscape) Color is used for the text in Inkscape, but the package 'color.sty' is not loaded}%
    \renewcommand\color[2][]{}%
  }%
  \providecommand\transparent[1]{%
    \errmessage{(Inkscape) Transparency is used (non-zero) for the text in Inkscape, but the package 'transparent.sty' is not loaded}%
    \renewcommand\transparent[1]{}%
  }%
  \providecommand\rotatebox[2]{#2}%
  \newcommand*\fsize{\dimexpr\f@size pt\relax}%
  \newcommand*\lineheight[1]{\fontsize{\fsize}{#1\fsize}\selectfont}%
  \ifx\svgwidth\undefined%
    \setlength{\unitlength}{276.52584434bp}%
    \ifx\svgscale\undefined%
      \relax%
    \else%
      \setlength{\unitlength}{\unitlength * \real{\svgscale}}%
    \fi%
  \else%
    \setlength{\unitlength}{\svgwidth}%
  \fi%
  \global\let\svgwidth\undefined%
  \global\let\svgscale\undefined%
  \makeatother%
  \begin{picture}(1,0.43808309)%
    \lineheight{1}%
    \setlength\tabcolsep{0pt}%
    \put(0,0){\includegraphics[width=\unitlength,page=1]{counterpoint.pdf}}%
  \end{picture}%
\endgroup%

%% file: musicology_name_1.bbl
\providecommand{\bysame}{\leavevmode\hbox to3em{\hrulefill}\thinspace}
\providecommand{\MR}{\relax\ifhmode\unskip\space\fi MR }
\providecommand{\MRhref}[2]{%
  \href{http://www.ams.org/mathscinet-getitem?mr=#1}{#2}
}
\providecommand{\href}[2]{#2}
\begin{thebibliography}{10}

\bibitem{OAAAthesis}
Octavio~A. Agust\'{i}n-Aquino, \emph{Extensiones microtonales de contrapunto
  (english version)}, Ph.D. thesis, Universidad Nacional Aut\'{o}noma de
  M\'{e}xico, 2011.

\bibitem{octaviotod}
Octavio~A. Agust\'{\i}n-Aquino, \emph{Tod der {H}omophonie! {E}s lebe der
  {K}ontrapunktsatz!}, \url{https://www.youtube.com/watch?v=2nh8xoLwC0o},
  November 2015.

\bibitem{Octavio}
Octavio~A. Agust{\'\i}n-Aquino, Julien Junod, and Guerino Mazzola,
  \emph{Computational counterpoint worlds: Mathematical theory, software, and
  experiments}, Computational Music Science, Springer International Publishing,
  Cham, Switzerland, 2015.

\bibitem{Mazzpolem}
Octavio~A. Agust{\'i}n-Aquino and Guerino Mazzola, \emph{{On D. Tymoczko’s
  Critique of Mazzola’s Counterpoint Theory}}, Memoirs of the Fourth
  International Seminar on Mathematical Music Theory (Octavio~A.
  Agust{\'i}n-Aquino and Emilio Lluis-Puebla, eds.), Sociedad Matemática
  Mexicana, México, 2011, pp.~43--48.

\bibitem{secsp}
Octavio~A. Agust{\'i}n-Aquino and Guerino Mazzola, \emph{A projection-oriented
  mathematical model for second-species counterpoint}, Preprint (2018).

\bibitem{Scriabinworld}
\bysame, \emph{{Contrapuntal Aspects of the Mystic Chord and Scriabin’s Piano
  Sonata No. 5}}, Mathematics and Computation in Music (Mariana Montiel,
  Francisco Gomez-Martin, and Octavio~A. Agust{\'i}n-Aquino, eds.), Springer
  International Publishing, Cham, Switzerland, 2019, pp.~3--20.

\bibitem{Theo}
Juan~S. Arias-Valero, Octavio~A. Agust{\'i}n-Aquino, and Emilio Lluis-Puebla,
  \emph{On first-species counterpoint theory}, Preprint (2021).

\bibitem{Fux}
Johann~J. Fux, Alfred Mann, and John Edmunds, \emph{The study of counterpoint
  from johann joseph fux’s gradus ad parnassum}, W. W. Norton, New York,
  1965.

\bibitem{Hartshorne}
Robin Hartshorne, \emph{Algebraic geometry}, Graduate Texts in Mathematics,
  Springer-Verlag, New York, 1977.

\bibitem{Jep}
Knud Jeppesen, \emph{Counterpoint: The polyphonic vocal style of the sixteenth
  century}, Dover Publications, New York, 1992.

\bibitem{Junod}
Julien Junod, \emph{Counterpoint worlds and morphisms : A graph-theoretical
  approach and its implementation on the rubato composer software}, Ph.D.
  thesis, University of Zurich, 2010.

\bibitem{Beau}
Guerino Mazzola, \emph{La v\' erit\' e du beau dans la musique}, Editions
  Delatour France, Paris, 2007.

\bibitem{MazzolaTopos}
Guerino Mazzola et~al., \emph{The topos of music: Geometric logic of concepts,
  theory, and performance}, Birkh\"auser Verlag, Basel, 2002. \MR{1938949
  (2004a:00013)}

\bibitem{MazzMuzz}
Guerino Mazzola and Daniel Muzzulini, \emph{{Deduktion des
  Quintparallelenverbots aus der Konsonanz-Dissonanz-Dichotomie}}, Musiktheorie
  (accepted) (1990).

\bibitem{Inicount}
Guerino Mazzola, Heinz-Gregor Wieser, Vreni Brunner, and Daniel Muzzulini,
  \emph{{A Symmetry-Oriented Mathematical Model of Classical Counterpoint and
  Related Neurophysiological Investigations by Depth EEG}}, Computers \&
  Mathematics with Applications \textbf{17} (1989), no.~4--6, 539--594.

\bibitem{Nieto}
Alejandro Nieto, \emph{Una aplicaci\'{o}n del teorema del contrapunto},
  Bachelor's thesis, ITAM, 2010.

\bibitem{Tittel}
Ernst Tittel, \emph{Der neue gradus}, Doblinger, Wien-M{\"u}nchen, 1959.

\bibitem{Tympolem}
Dmitri Tymoczko, \emph{{Mazzola’s Model of Fuxian Counterpoint}}, Mathematics
  and Computation in Music (Carlos Agon, Moreno Andreatta, G{\'e}rard Assayag,
  Emmanuel Amiot, Jean Bresson, and John Mandereau, eds.), Lecture Notes in
  Computer Science, vol. 6726, Springer, Heidelberg, 2011, pp.~297--310.
  \MR{2830991 (2012h:00024)}

\end{thebibliography}
